\documentclass[12pt, reqno]{amsart}
	\usepackage{amssymb, latexsym, amsmath, amsfonts}
	\usepackage{bm}
	\usepackage{mathrsfs}
	\usepackage[pdftex]{graphicx}
	\usepackage{float}
	\usepackage[numbers]{natbib}
	\usepackage{tikz}
	\usepackage{hyperref}
	\usepackage{changepage}
	
\newtheorem{thm}{Theorem}[section]

\newtheorem{cor}[thm]{Corollary}
\theoremstyle{definition}
\newtheorem{defn}[thm]{Definition}
\theoremstyle{remark}
\newtheorem{rem}[thm]{Remark}
\newtheorem{ex}{Example}

\newcommand{\mbb}{\mathbb}

\newcommand{\ov}{\overline}

\newcommand{\Om}{\Omega}

\newcommand{\ti}{\tilde}

\hypersetup{
	colorlinks,
	citecolor=blue,
	filecolor=black,
	linkcolor=red,
	urlcolor=black
}

\newcommand{\C}{\mathbb{C}} 
\newcommand{\R}{\mathbb{R}}
\newcommand{\Z}{\mathbb{Z}}

\newcommand{\sm}{\setminus}

\begin{document}
\title{A note on kernel functions of Dirichlet spaces}
	
\author{Sahil Gehlawat*,  Aakanksha Jain$^{\dagger}$ and Amar Deep Sarkar}

\address{SG: Universit\'{e} de Lille, Laboratoire de Math\'{e}matiques Paul Painlev\'{e}, CNRS U.M.R. 8524, 59655
Villeneuve d’Ascq Cedex, France.}
\email{sahil.gehlawat@univ-lille.fr, sahil.gehlawat@gmail.com}
	
\address{AJ: Department of Mathematics, Indian Institute of Science, Bangalore 560 012, India}
\email{aakankshaj@iisc.ac.in}

\address{ADS: School of Basic Sciences, Indian Institute of Technology Bhubaneswar, Argul 752050, India}
\email{amar@iitbbs.ac.in}

	

\thanks{*The author is supported by the Labex CEMPI (ANR-11-LABX-0007-01).}
\thanks{$\dagger$ The author is supported by the PMRF Ph.D. fellowship of the Ministry of Education, Government of India.}
\begin{abstract}
For a planar domain $\Om$, we consider the Dirichlet spaces with respect to a base point $\zeta\in\Om$ and the corresponding kernel functions. It is not known how these kernel functions behave as we vary the base point. In this note, we prove that these kernel functions vary smoothly. As an application of the smoothness result, we prove a Ramadanov-type theorem for these kernel functions on $\Om\times\Om$. This extends the previously known convergence results of these kernel functions. In fact, we have made these observations in a more general setting, that is, for weighted kernel functions and their higher-order counterparts.
\end{abstract}

\maketitle

\section{Introduction}

For a planar domain $\Om\subset\mbb {C}$ and $\zeta\in \Om$, the collection of holomorphic functions vanishing at $\zeta$ with $L^2$-integrable derivatives is called the Dirichlet space based at point $\zeta$. One can associate this space with the space of $L^2$- integrable holomorphic functions on $\Om$ which admits a primitive -- called the reduced Bergman space (see Definition 1.3 in \cite{GJS2}).

\medskip

Prior to this note, we proved a transformation formula for the weighted reduced Bergman kernels under proper holomorphic maps between bounded planar domains and also saw some applications of the transformation formula (see \cite{GJS1}). Subsequently, we considered the higher-order counterparts of the reduced Bergman kernel and studied various important properties of the same (see \cite{GJS2}). More specifically, we proved some Ramadanov-type theorems for these higher-order reduced Bergman kernels, and made significant observations about the boundary behaviour of these kernels. This note is a continuation of our previous efforts in studying the reduced Bergman kernel and related objects.

\medskip
The aim of this note is to study the $n$-th order weighted kernel functions $M_{\Om, \mu, n}(\cdot, \cdot)$ with weight $\mu$ associated with the Dirichlet space (see Definition \ref{Higher Order kernels}). These kernel functions are used to define the $n$-th order weighted reduced Bergman kernels. In 1978, M. Sakai (see \cite{MS1}) used these kernel functions on a Riemann surface $R$ in order to prove some fundamental results on the dimension of $AD(R)$ (the space of holomorphic functions on $R$ with finite Dirichlet integrals). In order to achieve this, he proved that for each $w \in R$, the sup norm of the kernel function $M_{R,n}(\cdot,w)$ is bounded above by $\frac{d M_{R,n}}{dz}(w,w)$, which in turn tells us that $M_{R,n}(\cdot,w)$ is bounded for each $w \in R$. Further works on these kernel functions and the span metric (the metric induced by $M_{R,n}(\cdot,\cdot)$) along similar lines can be found in \cite{Burbea1}, \cite{Burbea2} and \cite{MS2}. In this note, we study the regularity of these kernel functions. From the Definition \ref{Higher Order kernels}, it is clear that $M_{\Om,\mu,n}(z,w)$ is holomorphic as a function of $z \in \Om$. But the regularity corresponding to the other variable $w \in \Om$ is not a priori clear. We will prove that these kernel functions $M_{\Om,\mu,n}(\cdot, \cdot)$ are smooth in the complement of a very small set in $\Om \times \Om$. We will see by an example that these kernel functions need not be anti-holomorphic with respect to the second variable. Using the above regularity of these weighted kernel functions, we will prove a Ramadanov-type theorem for the $n$-th order weighted kernel functions given that it holds for the weighted reduced Bergman kernel. As a special case, we obtain a Ramadanov-type theorem for these kernel functions for eventually increasing sequence of domains. This extends observations made in Proposition 5.1 in \cite{MS2}, and Corollary 1.8 in \cite{GJS2}. One shall note that more substantial observations can be made for the $1$-st order weighted kernel functions $M_{\Om, \mu}$ as compared to the higher order kernel functions $M_{\Om,\mu,n}$ for $n >1$ -- similar to the case for the weighted reduced Bergman kernels in \cite{GJS2}.

\medskip
\noindent Before we define these weighted kernel functions, we shall see the type of weights that we will be working with throughout this article.

\begin{defn} (Z. Pasternak-Winiarski, \cite{PW1, PW2})
Let $\Omega\subset\mathbb{C}$ be a domain and $\mu$ be a positive measurable real-valued function on $\Omega$. The weight $\mu$ is called an admissible weight on $\Omega$ if for every compact set $K\subset\Omega$, there exists a constant $C_K > 0$ such that
\[
\sup_{z\in K} \vert f(z)\vert \leq C_K \Vert f\Vert_{L^2_{\mu}(\Omega)}
\]
for all $f\in\mathcal{O}(\Omega)\cap L^2_{\mu}(\Omega)$.
The space of admissible weights on $\Omega$ is denoted by $AW(\Omega)$. 
\end{defn}

\begin{rem}
It is known that if $\mu^{-a}$ is locally integrable on $\Om$ for some $a>0$, then $\mu\in AW(\Om)$.   
\end{rem}

\begin{defn}\label{Higher Order kernels} (See \cite{Bergman, MS1})
Let $\Om\subset\mathbb{C}$ be a domain, $\mu\in AW(\Om)$, $\zeta\in \Om$, and $n$ be a positive integer. The $n$-th order Dirichlet space based at point $\zeta$ is defined as
\[
AD^{\mu}(\Om,\zeta^n)=\left\lbrace f\in \mathcal{O}(\Om): f(\zeta)=f'(\zeta)=\cdots = f^{(n-1)}(\zeta)=0,\,\int_{\Om}\lvert f'(z)\rvert^2 \mu(z) dA(z) <\infty\right\rbrace.
\]
This is a Hilbert space with respect to the inner product
\[
\langle f,g\rangle_{AD^{\mu}(\Om,\zeta^n)}=\int_{\Om}f'(z)\,\overline{g'(z)}\mu(z)\,dA(z),\quad f,g\in AD^{\mu}(\Om,\zeta^n).
\]
The linear functional defined by
$
AD^{\mu}(\Om,\zeta^n)\ni f\mapsto f^{(n)}(\zeta)\in\mathbb{C},
$
is continuous. By Riesz representation theorem, there exists a unique function $M_{\Om,\mu,n}(\cdot,\zeta)\in AD^{\mu}(\Om,\zeta^n)$ such that
$f^{(n)}(\zeta)=\langle f, M_{\Om,\mu,n}(\cdot,\zeta)\rangle_{AD^{\mu}(\Om,\zeta^n)}$ for every $f\in AD^{\mu}(\Om,\zeta^n)$.

\medskip

The function $M_{\Om,\mu,n}(\cdot,\cdot)$ is called the $n$-th order weighted kernel function associated with the Dirichlet space with respect to weight $\mu$. Define
\[
\tilde{K}_{\Om,\mu,n}(z,\zeta)=\frac{\partial}{\partial z} M_{\Om,\mu,n}(z,\zeta),\quad z,\zeta\in \Om.
\]
The kernel $\tilde{K}_{\Om,\mu,n}(\cdot,\cdot)$ is called the $n$-th order weighted reduced Bergman kernel of $\Om$ with respect to the weight $\mu$. Putting $n=1$ gives the weighted reduced Bergman kernel $\tilde{K}_{\Om,\mu}(\cdot,\cdot)$ of $\Om$ with weight $\mu$. 
\end{defn}

\begin{rem}
The weighted reduced Bergman kernel $\tilde{K}_{\Om,\mu}(\cdot,\cdot)$ is the reproducing kernel of a closed subspace of the weighted Bergman space. This is the space of all $L^2(\mu)$-integrable holomorphic functions on $\Om$ which admits a primitive. Therefore, $\tilde{K}_{\Om,\mu}(\cdot,\cdot)$ is holomorphic in the first variable and anti-holomorphic in the second variable. This in turn implies that $\tilde{K}_{\Om,\mu}\in C^{\infty}(\Om\times\Om)$.

\medskip

\noindent It is known (see \cite[p. 26]{Bergman}, \cite[p. 476]{Burbea2}) that for a domain $\Om\subset\mathbb{C}$, $\mu\in AW(\Om)$, and $n\geq 2$, 
\begin{equation}\label{determinant}
\tilde{K}_{\Om,\mu,n}(z,\zeta)=
\frac{(-1)^{n-1}}{J_{n-2}}
\det\left(
\begin{matrix}
\tilde{K}_{0,\bar{0}}(z,\zeta)& \ldots & \tilde{K}_{0,\overline{n-1}}(z,\zeta)\\
\tilde{K}_{0,\bar{0}}& \ldots & \tilde{K}_{0,\overline{n-1}}\\
\tilde{K}_{1,\bar{0}}& \ldots & \tilde{K}_{1,\overline{n-1}}\\
\vdots & &\vdots\\
\tilde{K}_{n-2,\bar{0}}& \ldots & \tilde{K}_{n-2,\overline{n-1}}
\end{matrix}
\right),
\end{equation}
where $J_{n}=\det\left(\tilde{K}_{j\bar{k}}\right)_{j,k=0}^n$ and 
\[
\tilde{K}_{j\bar{k}}(z,\zeta)=\frac{\partial^{j+k}}{\partial z^j\partial\bar{\zeta}^k}\tilde{K}_{\Om, \mu}(z,\zeta),\quad \tilde{K}_{j\bar{k}}\equiv \tilde{K}_{j\bar{k}}(\zeta,\zeta).
\]
Here $J_n>0$ for all $\zeta\notin N_{\Om}(\mu):=\{z\in \Om: \tilde{K}_{\Om,\mu}(z,z)=0\}$. Thus, $\tilde{K}_{\Om,\mu,n}\in C^{\infty}(\Om\times (\Om\setminus N_{\Om}(\mu))$. As a special case, if $\mu\in L^1(\Om)$, then $N_{\Om}(\mu)=\emptyset$, and therefore $\tilde{K}_{\Om,\mu,n}\in C^{\infty}(\Om\times \Om)$.

\end{rem}

\begin{thm}\label{smoothness}
Let $ \Omega \subset \mathbb{C} $ be a domain, $\mu\in AW(\Omega)$ and $n$ be a positive integer. The kernel function $M_{\Omega,\mu}(\cdot, \cdot) \in C^{\infty}(\Om\times\Om)$. For $n>1$, the higher-order kernel function $M_{\Omega,\mu,n}(\cdot, \cdot) \in C^{\infty}(\Om\times(\Om\setminus N_{\Om}(\mu)))$, where $N_{\Omega}(\mu) = \{\zeta\in\Om : \tilde{K}_{\Om,\mu}(\zeta,\zeta) = 0\}$.

\smallskip

\noindent Moreover, for $n\geq 1$ and non-negative integers $r,s$
\[
\frac{\partial^{r+s} M_{\Om,\mu,n}(z,\zeta)}{\partial\zeta ^r\partial{\ov{\zeta}}^s} 
=
\begin{cases}
\int\limits_{\zeta}^{z}\frac{\partial^{s} \tilde K_{\Omega,\mu,n}(\xi, \zeta)}{\partial{\ov{\zeta}}^s} d\xi 
& \text{for } r=0 \\
\int\limits_{\zeta}^{z}\frac{\partial^{r+s} \tilde K_{\Omega,\mu,n}(\xi, \zeta)}{\partial \zeta^r\partial{\ov{\zeta}}^s} d\xi 
-
\sum\limits_{k =0}^{r-1}\frac{\partial^{k+s} \tilde K_{\Om,\mu,n}^{(\overline{r-1-k})}}{\partial \zeta^k\partial{\ov{\zeta}}^s} 
& \text{for } r\geq 1
\end{cases}
\]
where for a positive integer $m$
\[
\tilde K_{\Om,\mu,n}^{(\overline{m})}(z, \zeta) = \frac{\partial^m \tilde K_{\Om,\mu,n}(z, \zeta)}{\partial \zeta^m}
\quad
\text{and}
\quad
\tilde K_{\Om,\mu,n}^{(\overline{m})} = \tilde K_{\Om,\mu,n}^{(\overline{m})}(\zeta, \zeta).
\]
\end{thm}

\medskip
\noindent The following example gives the expression for the $n$-th order kernel function for the unit disc $\mbb{D}$.

\begin{ex}
Let $\zeta \in \mbb{D}$ and $f \in AD(\mbb{D}, \zeta^n)$ for $n \in \mbb{Z}^+$, that is $f^{(k)}(\zeta) = 0$ for all $0 \le k \le n-1$. For $g\in AD(\mbb{D},0^n)$ and $0<r<1$, the Cauchy integral formula gives us 

\begin{eqnarray*}
g^{(n)}(0)
&=& \frac{(n-1)!}{2\pi i} \int_{\vert \xi\vert = r} \frac{g'(\xi)}{\xi^n} \, d\xi = \frac{(n-1)!}{2\pi i} \int_{0}^{2\pi} \frac{g'(re^{it})}{(re^{it})^n} r e^{it} i \, dt
\\
&=&
\frac{(n-1)!}{2\pi} \int_{0}^{2\pi} \frac{g'(re^{it})}{r^{n-1}} \overline{(e^{it})^{n-1}} \, dt.
\end{eqnarray*}

\noindent Multiplying both sides by $r^{2n-1}$ and integrating with respect to parameter $r$, we get
\[
\int_{0}^{1} g^{(n)}(0) r^{2n-1} dr = \frac{(n-1)!}{2\pi} \int_{0}^{1} \int_{0}^{2\pi} g'(re^{it}) \overline{(re^{it})^{n-1}} \,r \, dr \, dt.
\]
By change of variables on the right hand side of the above equation, we get

\begin{equation}\label{Eqn3.1}
g^{(n)}(0) = \frac{n!}{\pi} \int_{\mbb{D}} g'(\xi) \overline{(\xi)^{n-1}} \, dA(\xi).
\end{equation}

\noindent Let $\tilde{K}_{n}(\cdot, \cdot)$ denote the $n$-th order reduced Bergman kernel of $\mbb{D}$. For all $\zeta \in \mbb{D}$, and $f \in AD(\mbb{D}, \zeta^n)$, we have
\begin{equation}\label{Eqn3.2}
f^{(n)}(\zeta) = \int_{\mbb{D}} f'(\xi) \overline{\tilde{K}_{n}(\xi,\zeta)} \, dA(\xi).
\end{equation}

\noindent Let $\phi_{\zeta}: \mbb{D} \rightarrow \mbb{D}$ be the automorphism of unit disc given by $\phi_{\zeta}(z) = \frac{\zeta - z}{1- z\overline{\zeta}}$. Note that $\phi_{\zeta}(0) = \zeta$, $\phi_{\zeta}(\zeta) = 0$, and $\phi_{\zeta} \circ \phi_{\zeta}(z) = z$ for all $z \in \mbb{D}$. Consider the holomorphic function $f \circ \phi_{\zeta}$. Observe that $f \circ \phi_{\zeta} \in AD(\mbb{D}, 0^n)$ and $(f \circ \phi_{\zeta})^{(n)}(0) = f^{(n)}(\zeta) (\phi_{\zeta}'(0))^n$. On substituting $g=f\circ \phi_{\zeta}$ in equation (\ref{Eqn3.1}), we get

\begin{eqnarray*}
f^{(n)}(\zeta) (\phi_{\zeta}'(0))^n &=& \frac{n!}{\pi} \int_{\mbb{D}} (f\circ \phi_{\zeta})'(\chi) \overline{(\chi)^{n-1}} \, dA(\chi) 
\\
&=& \frac{n!}{\pi} \int_{\mbb{D}} f'(\phi_{\zeta}(\chi)) \phi_{\zeta}'(\chi) \overline{(\chi)^{n-1}} \, dA(\chi). \\
\end{eqnarray*}

\noindent Now by doing change of variables $\xi = \phi_{\zeta}(\chi)$, we get $\chi = \phi_{\zeta}(\xi)$, and $\phi_{\zeta}'(\chi) = \frac{1}{\phi_{\zeta}'(\xi)}$. Therefore

\begin{eqnarray*}
f^{(n)}(\zeta) (\phi_{\zeta}'(0))^n &=& \frac{n!}{\pi} \int_{\mbb{D}} f'(\xi) (\phi_{\zeta}'(\xi))^{-1} \overline{(\phi_{\zeta}(\xi))^{n-1}} \, \vert \phi_{\zeta}'(\xi)\vert^2 \, dA(\xi) \\
&=& \frac{n!}{\pi} \int_{\mbb{D}} f'(\xi) \overline{(\phi_{\zeta}(\xi))^{n-1}} \, \overline{\phi_{\zeta}'(\xi)} \, dA(\xi). \\
\end{eqnarray*}
Therefore, we get

\begin{equation}\label{Eqn3.3}
f^{(n)}(\zeta) = \int_{\mbb{D}} f'(\xi) \frac{n! \overline{(\phi_{\zeta}(\xi))^{n-1}} \, \overline{\phi_{\zeta}'(\xi)}}{\pi (\phi_{\zeta}'(0))^{n}} \, dA(\xi).
\end{equation}

\noindent Comparing the equations (\ref{Eqn3.2}) and (\ref{Eqn3.3}), we get 
\[
\tilde{K}_{n}(\xi, \zeta) = \frac{n!}{\pi} \frac{(\phi_{\zeta}(\xi))^{n-1} \, \phi_{\zeta}'(\xi)}{\overline{(\phi_{\zeta}'(0))^{n}}}.
\]

\noindent Since $\frac{d M_{n}(\xi, \zeta)}{d \xi} = \tilde{K}_{n}(\xi, \zeta)$ and $M_{n}(\zeta, \zeta) = 0$, therefore
\[
M_{n}(\xi, \zeta) = \frac{(n-1)!}{\pi} \frac{(\phi_{\zeta}(\xi))^{n}}{\overline{(\phi_{\zeta}'(0))^{n}}}.
\]

\noindent We can check that $\phi_{\zeta}'(\xi) = \frac{\vert \zeta\vert^2 - 1}{(1 - \overline{\zeta}\xi)^2}$, which gives $\phi_{\zeta}'(0) = \vert \zeta\vert^2 - 1$. Therefore,
\[
M_{n}(\xi, \zeta) = \frac{(n-1)!}{\pi} \frac{(\zeta - \xi)^{n}}{(1 - \overline{\zeta}\xi)^{n}} \frac{1}{(\vert \zeta\vert^2 - 1)^n}.
\]
Thus, for $n \ge 1$, the $n$-th order kernel function is given by
\begin{equation}\label{RBKforUnitDisc}
M_{n}(\xi, \zeta) = \frac{(n-1)!}{\pi} \frac{(\xi - \zeta)^{n}}{(1 - \overline{\zeta}\xi)^{n} (1 - \vert \zeta\vert^2)^{n}}.
\end{equation}
\end{ex}

\medskip

Now suppose $\{\Om_j\}_{j \ge 1}$ be a sequence of planar domains with $\mu_{j} \in AW(\Om_j)$ and $n$ be a positive integer. Ramadanov \cite{Rama} showed that if $\Om_j \subset \Om_{j+1}$ for all $j \in \mathbb{Z}^+$, and $\Om := \bigcup_{j =1}^{\infty} \Om_j$, the Bergman kernel $K_{j}(\cdot, \cdot)$ corresponding to $\Om_j$ converges uniformly on compacts of $\Om \times \Om$ to the Bergman kernel $K(\cdot, \cdot)$ corresponding to $\Om$. The question here is to study the variation of the kernel functions $M_{\Om_j, \mu_j, n}$, given some type of convergence of the domains $(\Om_j, \mu_j)$.

\medskip
In 1979, M. Sakai proved that if $\{\Om_j\}_{j \ge1}$ is an increasing sequence of planar domains and $\mu_j \le \mu_{j+1}$ for all $j \ge 1$ where $\mu_{j} \in AW(\Om_j)$, and there exist a domain $\Om \subset \mathbb{C}$ with $\mu \in AW(\Om)$, such that $\Om = \cup_{j \ge 1} \Om_j$ and $\mu_j \to \mu$ pointwise, then for each $\zeta \in \Om$, $M_{\Om_{j},\mu_{j},n}(\cdot, \zeta) \to M_{\Om,\mu,n}(\cdot, \zeta)$ uniformly on compacts of $\Om$ (see Proposition 5.1 in \cite{MS2}). In \cite{GJS2}, we proved similar observations for the case when $\{(\Om_j, \mu_j)\}_{j \ge 1}$ is eventually increasing, that is, for each $j \in \mathbb{Z}^+$, there exist $k(j) \in \mathbb{Z}^+$ such that $\Om_j \subset \Om_{l}$ and $\mu_j \le \mu_{l}$ for all $l \ge k(j)$, with $\Om = \cup_{j \ge 1} \Om_j$ and $\mu_j \to \mu$. Here we want to talk about the variation of $M_{\Om_j, \mu_j, n}(\cdot,\cdot)$ on $\Om \times \Om$. In fact, Corollary \ref{C:Ramadanov} tells us that the convergence is uniform on compacts of the complement of a thin set in $\Om \times \Om$.

\begin{thm}\label{Ramadanov}
Let $\Om,\, \Om_j\subset \mathbb{C}$ be domains with $\mu_{j} \in AW(\Om_j)$, $\mu \in AW(\Om)$ and $n$ be a positive integer. Assume that every compact set $K\subset\Omega$ is eventually contained in $\Om_j$. If 
\[
\lim_{j\rightarrow\infty}  \tilde K_{\Om_j,\mu_j}(z, \zeta) = \tilde K_{\Om,\mu}(z, \zeta) 
\]
locally uniformly on $ \Om \times \Om $, then 
\begin{enumerate}
\item[1.] the sequence of kernel functions $ M_{\Om_j,\mu_j}(\cdot, \cdot) $ converges to the kernel function $ M_{\Om,\mu}(\cdot, \cdot) $ locally uniformly on $ \Om \times \Om $. Moreover, all partial derivatives of $ M_{\Om_j,\mu_j} $ converge to the corresponding partial derivatives of $ M_{\Om,\mu} $ locally uniformly on $ \Om \times \Om$.
\item[2.] for $n>1$, the sequence of $n$-th order kernel functions $ M_{\Om_j,\mu_j,n}(\cdot, \cdot) $ converges to the $n$-th order kernel function $ M_{\Om,\mu,n}(\cdot, \cdot) $ locally uniformly on $ \Om \times (\Om\sm N_{\Om}(\mu)) $. Moreover, all partial derivatives of $ M_{\Om_j,\mu_j,n} $ converge to the corresponding partial derivatives of $ M_{\Om,\mu,n} $ locally uniformly on $ \Om \times (\Om\sm N_{\Om}(\mu))$.
\end{enumerate}
\end{thm}

\medskip
\noindent Now suppose that $\{(\Om_j, \mu_j)\}_{j \ge 1}$ is an eventually increasing sequence of planar domains, and there exist a planar domain $\Om$ with $\mu \in AW(\Om)$ such that $\Om = \cup_{j \ge 1} \Om_j$ and $\mu_j \to \mu$. Using Theorem \ref{Ramadanov} and Theorem 1.6 in \cite{GJS2}, the following corollary is immediate.

\begin{cor}\label{C:Ramadanov}
	Suppose that the sequence of domains $ \Om_j$ increases eventually to $ \Om $ and $ \mu_j $ increases eventually to $ \mu $ as $ j \to \infty $. Then the sequence of kernel functions $ M_{\Om_j,\mu_j,n} $ of the domain $ \Om_j $ converges to the kernel function $ M_{\Om,\mu,n} $ uniformly on compact subsets of $ \Om \times (\Om\sm N_{\Om}(\mu)) $.
\end{cor}


\section{Proof of Theorem \ref{smoothness}}

\begin{proof}
Let $n>1$. For $z,w\in\Omega$, let $\int_z^w$ denote the integration along a path from $z$ to $w$ in $\Omega$. Recall that
\[
\frac{\partial}{\partial z}M_{\Omega,\mu,n}(z,\zeta) = \tilde{K}_{\Omega,\mu,n}(z,\zeta)
,\quad\text{and}\quad
M_{\Omega,\mu,n}(\zeta, \zeta) = 0\quad \text{ for }z,\zeta\in\Omega.
\]
Therefore,
\[
M_{\Om,\mu,n}(z, \zeta) = \int\limits_{\zeta}^z\tilde K_{\Om,\mu,n}(\xi, \zeta) d\xi,\quad z,\zeta\in\Om.
\]
Since $M_{\Om,\mu,n}$ is a primitive of $\tilde{K}_{\Omega,\mu,n}$, the above integral does not depend upon the choice of path. Note that $\tilde{K}_{\Om,\mu,n}(\cdot,\cdot)\in C^{\infty}(\Om\times(\Om\setminus N_{\Om}(\mu)))$. Fix $z\in\Omega$ and let $\zeta\in \Om\setminus N_{\Om}(\mu)$. For $ w \in \C $ with small enough modulus,
\begin{multline*}
M_{\Om,\mu,n}(z,\zeta + w) - M_{\Om,\mu,n}(z,\zeta) 
=
\int_{\zeta + w}^{z}\tilde K_{\Omega,\mu,n}(\xi, \zeta + w) d\xi -\int_{\zeta}^{z}\tilde K_{\Omega,\mu,n}(\xi, \zeta) d\xi
\\=
\int_{\zeta}^{z}\left(\tilde K_{\Omega,\mu,n}(\xi, \zeta + w) - \tilde K_{\Om,\mu,n}(\xi, \zeta)\right) d\xi 
-
\int_{\zeta}^{\zeta + w}\tilde K_{\Omega,\mu,n}(\xi, \zeta + w) d\xi.
\end{multline*}
Let $\zeta = u + iv$. Take $w=h$ with $h\in\mathbb{R}$. Since $\tilde{K}_{\Om,\mu,n}(\cdot,\cdot)\in C^{\infty}(\Om\times(\Om\setminus N_{\Om}(\mu)))$, an application of DCT gives
\begin{eqnarray*}
\lim_{h \to 0}\int\limits_{\zeta}^{z}\left(\frac{\tilde K_{\Omega,\mu,n}(\xi, \zeta + h) - \tilde K_{\Omega,\mu,n}(\xi, \zeta)}{h}\right) d\xi
&=&
\int\limits_{\zeta}^{z}\lim_{h \to 0}\left(\frac{\tilde K_{\Omega,\mu,n}(\xi, \zeta + h) - \tilde K_{\Omega,\mu,n}(\xi, \zeta)}{h}\right) d\xi
\\&=&
\int\limits_{\zeta}^{z}\frac{\partial \tilde K_{\Omega,\mu,n}(\xi, \zeta)}{\partial u} d\xi.
\end{eqnarray*}
By taking the curve $ \gamma(t) = \zeta + th $, $ t \in [0, 1] $ (for small enough $ h $), we get
\[
\lim_{h \to 0} \frac{1}{h}\left(\int\limits_{\zeta}^{\zeta + h}\tilde K_{\Omega,\mu,n}(\xi, \zeta + h) d\xi\right) = \lim_{h\rightarrow 0}\int\limits_{0}^{1} \tilde K_{\Omega,\mu,n}(\zeta + th, \zeta + h) dt = \tilde K_{\Omega,\mu,n}(\zeta, \zeta).
\]
The last equality follows from the continuity of $ \tilde K_{\Om,\mu,n}(z, \zeta) $ in both the variables on $\Om \times (\Om\setminus N_{\Om})$.
Similarly, for $w= ih$ with $ h \in \R $, we obtain	
\begin{eqnarray*}
\lim_{h \to 0}\int\limits_{\zeta}^{z}\left(\frac{\tilde K_{\Omega,\mu,n}(\xi, \zeta + ih) - \tilde K_{\Omega,\mu,n}(\xi, \zeta)}{h}\right) d\xi
&=&
\int\limits_{\zeta}^{z}\lim_{h \to 0}\left(\frac{\tilde K_{\Omega,\mu,n}(\xi, \zeta + ih) - \tilde K_{\Omega,\mu,n}(\xi, \zeta)}{h}\right) d\xi
\\&=&
\int\limits_{\zeta}^{z}\frac{\partial \tilde K_{\Omega,\mu,n}(\xi, \zeta)}{\partial v} d\xi,
\end{eqnarray*}
and taking $ \gamma(t) = \zeta +ith $ for $ t \in [0, 1] $, we have
\[
\lim_{h \to 0} \frac{1}{h}\left(\int\limits_{\zeta}^{\zeta + ih}\tilde K_{\Omega,\mu,n}(\xi, \zeta + ih) d\xi\right) = i \lim_{h \to 0} \int\limits_{0}^{1} \tilde K_{\Omega,\mu,n}(\zeta + ith, \zeta + h) dt = i\tilde K_{\Omega,\mu,n}(\zeta, \zeta).
\]
Thus, we obtain from above calculations that
\begin{eqnarray*}
\frac{\partial M_{\Om,\mu,n}(z,\zeta)}{\partial \zeta}
&=&
\int\limits_{\zeta}^{z}\frac{\partial \tilde K_{\Omega,\mu,n}(\xi, \zeta)}{\partial \zeta} d\xi - \frac{1}{2}\left(\tilde K_{\Omega,\mu,n}(\zeta, \zeta) + \tilde K_{\Omega,\mu,n}(\zeta, \zeta)\right) 
\\&=&
\int\limits_{\zeta}^{z}\frac{\partial \tilde K_{\Omega,\mu,n}(\xi, \zeta)}{\partial \zeta} d\xi - \tilde K_{\Omega,\mu,n}(\zeta, \zeta).
\end{eqnarray*}
Also,
\begin{eqnarray*}
\frac{\partial M_{\Om,\mu,n}(z,\zeta)}{\partial \bar{\zeta}}
&=&
\int\limits_{\zeta}^{z}\frac{\partial \tilde K_{\Omega,\mu,n}(\xi, \zeta)}{\partial \bar{\zeta}} d\xi - \frac{1}{2}\left(\tilde K_{\Omega,\mu,n}(\zeta, \zeta) - \tilde K_{\Omega,\mu,n}(\zeta, \zeta)\right) 
\\&=&
\int\limits_{\zeta}^{z}\frac{\partial \tilde K_{\Omega,\mu,n}(\xi, \zeta)}{\partial \bar{\zeta}} d\xi.
\end{eqnarray*}
It now follows from induction and the fact that $\tilde{K}_{\Om,\mu,n}(\cdot,\cdot)\in C^{\infty}(\Om\times(\Om\setminus N_{\Om}(\mu)))$, for positive integers $r,s$, all the partial derivatives in $\zeta,\overline{\zeta}$ commutes, and 
\begin{eqnarray*}
\frac{\partial^{r+s} M_{\Om,\mu,n}(z,\zeta)}{\partial\zeta ^r\partial{\bar{\zeta}}^s} 
&=&
\frac{\partial^{r+s} M_{\Om,\mu,n}(z,\zeta)}{\partial{\bar{\zeta}}^s\partial\zeta ^r}
\\&=&
\frac{\partial^s}{\partial {\bar{\zeta}}^s} 
\left(
\int\limits_{\zeta}^{z}\frac{\partial^r \tilde K_{\Omega,\mu,n}(\xi, \zeta)}{\partial \zeta^r} d\xi 
-
\sum_{k =0}^{r-1}\frac{\partial^k \tilde K_{\Om,\mu,n}^{(\overline{r-1-k})}}{\partial \zeta^k}
\right)
\\&=&
\int\limits_{\zeta}^{z}\frac{\partial^{r+s} \tilde K_{\Omega,\mu,n}(\xi, \zeta)}{\partial{\bar{\zeta}}^s\partial \zeta^r} d\xi 
-
\sum_{k =0}^{r-1}\frac{\partial^{k+s} \tilde K_{\Om,\mu,n}^{(\overline{r-1-k})}}{\partial{\bar{\zeta}}^s\partial \zeta^k}.
\end{eqnarray*}
Moreover, the moment we differentiate with respect to $z$, the integrals disappear and the resulting expression is smooth by the smoothness property of $\tilde K_{\Om,\mu,n}$ on $\Om\times(\Om\setminus N_{\Om}(\mu))$.
Additionally, the functions involved are holomorphic in $z$. Thus, we have proved that $ M_{\Om,\mu, n}(z, \zeta) $ is a $ C^{\infty} $-smooth function on $ \Om \times (\Om\setminus N_{\Om}(\mu)) $.

\medskip

For $n=1$, it can be proved in a similar manner that the kernel function $M_{\Om,\mu}(z, \zeta)$ is smooth on $\Om \times \Om$ because $\tilde{K}_{\Om,\mu}(\cdot,\cdot)\in C^{\infty}(\Om\times\Om)$.
\end{proof}


\section{Proof of Theorem \ref{Ramadanov}}

\begin{proof}
Let $n>1$. Fix $ (z_0, \zeta_0) \in \Om \times (\Om\sm N_{\Om}(\mu)) $ and choose $ r> 0 $ such that $ \overline{B(z_0, r)} \times \overline{B(\zeta_0, r)}  \subset \Om \times (\Om\sm N_{\Om}(\mu))$. Let $ \gamma: [0, 1] \longrightarrow \Om \sm N_{\Om}(\mu) $ be a piecewise $ C^1 $-smooth curve such that $ \gamma(0) = \zeta_0 $ and  $ \gamma(1) = z_0 $. The set 
\[
W := (\gamma \cup \overline{B(z_0, r)} \cup \overline{B(\zeta_0, r)}) \times \overline{B(\zeta_0, r)}
\]
is a compact subset of $\Om \times (\Om \sm N_{\Om}(\mu))$. We may assume, without loss of generality, that $W \subset \Om_j \times (\Om_j\sm N_{\Om_j}(\mu_j)) $ for all $ j $. Now, for $ \zeta \in  B(\zeta_0, r)  $ and $ z \in  B(z_0, r)  $, define a path $ \sigma_{\zeta, z} = \gamma_{z}*\gamma* \gamma_{\zeta} $ joining $ \zeta $ and $ z $, where $ \gamma_{\zeta}(t) := \zeta + t(\zeta_0 - \zeta) $ and $ \gamma_{z}(t) := z_0 + t(z - z_0) $ for all $ t \in [0, 1] $. Set $l(\gamma):= length(\gamma)$. Observe that $ length(\sigma_{\zeta, z} )\leq 2r + l(\gamma)$ as
\[
\int_{\sigma_{\zeta, z}} |d\xi| \leq \int_{\sigma_{z}} |d\xi| + \int_{\gamma} |d\xi|  + \int_{\sigma_{\zeta}} |d\xi| \leq 2r + l(\gamma).
\]
By Theorem $\ref{smoothness}$, for non-negative integers $r,s$
\[
\frac{\partial^{r+s} M_{\Om_j,\mu_j,n}(z,\zeta)}{\partial\zeta ^r\partial{\ov{\zeta}}^s} 
=
\begin{cases}
\int\limits_{\zeta}^{z}\frac{\partial^{s} \tilde K_{\Omega_j,\mu_j,n}(\xi, \zeta)}{\partial{\ov{\zeta}}^s} d\xi 
& \text{for } r=0 \\
\int\limits_{\zeta}^{z}\frac{\partial^{r+s} \tilde K_{\Omega_j,\mu_j,n}(\xi, \zeta)}{\partial \zeta^r\partial{\ov{\zeta}}^s} d\xi 
-
\sum\limits_{k =0}^{r-1}\frac{\partial^{k+s} \tilde K_{\Om_j,\mu_j,n}^{(\overline{r-1-k})}}{\partial \zeta^k\partial{\ov{\zeta}}^s} 
& \text{for } r\geq 1.
\end{cases}
\]
By the determinant formula (\ref{determinant}), note that  the local uniform convergence of $ \tilde K_{\Om_j,\mu_j} $ to $ \ti K_{\Om,\mu} $ on $ \Om \times \Om $ implies that all the partial derivatives of $ \tilde K_{\Om_j,\mu_j,n} $ converges to the corresponding partial derivatives of $ \tilde K_{\Om,\mu,n} $  uniformly on compact subsets of $ \Om \times (\Om\sm N_{\Om}(\mu)) $. 

\smallskip

Let $r,s$ be fixed non-negative integers. Let $ \epsilon > 0 $. Since $ \frac{\partial^{r+s} \tilde K_{\Om_j,\mu_j,n}(\xi, \zeta)}{\partial \zeta^r\partial \ov \zeta^s} \to  \frac{\partial^{r+s} \tilde K_{\Om,\mu,n}(\xi, \zeta)}{\partial \zeta^r\partial \ov \zeta^s}  $ uniformly on compact subsets of $ \Om \times (\Om\sm N_{\Om}(\mu)) $, there exists $ j_0(r,s) \in \Z^+ $ such that 
\[
\sup_{(\xi, \zeta) \in W}\left\vert \frac{\partial^{r+s} \tilde K_{\Om_j,\mu_j,n}(\xi, \zeta)}{\partial \zeta^r\partial \ov \zeta^s} -  \frac{\partial^{r+s} \tilde K_{\Om,\mu,n}(\xi, \zeta)}{\partial \zeta^r\partial \ov \zeta^s} \right\vert < \frac{\epsilon}{2r + l(\gamma)}
\]
for all $ j \geq j_0(r,s) $. Therefore,
\begin{multline*}
\sup_{(z, \zeta) \in B(z_0, r) \times B(\zeta_0,r)}
\left\lvert
\int\limits_{\zeta}^{z}\frac{\partial^{r+s} \tilde K_{\Omega_j,\mu_j,n}(\xi, \zeta)}{\partial \zeta^r\partial{\ov{\zeta}}^s} d\xi
-
\int\limits_{\zeta}^{z}\frac{\partial^{r+s} \tilde K_{\Omega,\mu,n}(\xi, \zeta)}{\partial \zeta^r\partial{\ov{\zeta}}^s} d\xi
\right\rvert 
\\
=
\sup_{(z, \zeta) \in B(z_0, r) \times B(\zeta_0,r)} \left\vert\int_{\sigma_{\zeta, z}}\left(\frac{\partial^{r+s} \tilde K_{\Om_j,\mu_j,n}(\xi, \zeta)}{\partial \zeta^r \partial \ov \zeta^s} - \frac{\partial^{r+s} \tilde K_{\Om,\mu,n}(\xi, \zeta)}{\partial \zeta^r \partial \ov \zeta^s} \right)d\xi \right\vert \\
\le
\sup_{(\xi, \zeta) \in W}\left\vert \frac{\partial^{r+s} \tilde K_{\Om_j,\mu_j,n}(\xi, \zeta)}{\partial \zeta^r \partial \ov \zeta^s} -  \frac{\partial^{r+s} \tilde K_{\Om,\mu,n}(\xi, \zeta)}{\partial \zeta^r \partial \ov \zeta^s} \right\vert(2r + l(\gamma))
< \, \epsilon
\end{multline*}
for all $j\geq j_0(r,s)$.
Moreover, for all integers $r\geq 1$ and $s\geq 0$,
\[
\lim_{j\rightarrow\infty}
\left(
\sum\limits_{k =0}^{r-1}\frac{\partial^{k+s} \tilde K_{\Om_j,\mu_j,n}^{(\overline{r-1-k})}}{\partial \zeta^k\partial{\ov{\zeta}}^s} 
\right)
=
\sum\limits_{k =0}^{r-1}\frac{\partial^{k+s} \tilde K_{\Om,\mu,n}^{(\overline{r-1-k})}}{\partial \zeta^k\partial{\ov{\zeta}}^s} 
\]
uniformly for all $\zeta\in B(\zeta_0,r)$ as all the partial derivatives of $ \tilde K_{\Om_j,\mu_j,n} $ converges to the corresponding partial derivatives of $ \tilde K_{\Om,\mu,n} $  uniformly on compact subsets of $ \Om \times (\Om\sm N_{\Om}(\mu)) $. 

\medskip

As noted before, the moment we differentiate the kernel functions with respect to $z$, the integrals disappear and the resulting expression is a linear combination of partial derivatives of the n-th order reduced Bergman kernels. Additionally, all the kernel functions involved are holomorphic in $z$. 
Thus, we have proved that all the partial derivatives of $M_{\Om_j,\mu_j,n} $ converge to the corresponding partial derivatives of $ M_{\Om,\mu,n} $ locally uniformly on $ \Om \times (\Om\sm N_{\Om}(\mu))$.

\medskip

Exactly similar calculations and the fact that the reduced Bergman kernel is holomorphic in the first variable and anti-holomorphic in the second variable will lead us to conclude the local uniform convergence of all the partial derivatives of $M_{\Om_j,\mu_j} $ to the corresponding partial derivatives of $ M_{\Om,\mu} $ on $ \Om \times \Om$.	
\end{proof}




\end{document}